\documentclass{article}

\usepackage{amsmath,amsthm}
\usepackage{enumerate}
\usepackage{amssymb}
\usepackage{graphicx}
\usepackage{hyperref}
\usepackage[usenames,dvipsnames]{color}
\usepackage{indentfirst} 

\usepackage{tikz}

\newcommand{\degb}{\deg_{\textnormal{B}}}

\newcommand{\dom}{\textnormal{dom\,}}
\newcommand{\im}{\textnormal{Im\,}}
\newcommand{\coker}{\textnormal{coker\,}}
\newcommand{\R}{\mathbb{R}}

\numberwithin{equation}{section}

\theoremstyle{plain}
\newtheorem{theorem}{Theorem}[section]
\newtheorem{af}{Assertion}
\newtheorem{lema}{Lema}
\newtheorem{obs}{Remark}
\usepackage{amsfonts}
\author{Adriano Leando da Costa Peixoto\thanks{Email adress: alcp@ime.usp.br}\\ Instituto de Matemática e Estatística da USP}
\title{Existence of positive solution for a nonlinear problem with mixed conditions.}

\date{}

\usepackage[english]{babel} 
\usepackage{amssymb}

\theoremstyle{plain}

\begin{document}
	\maketitle	

\begin{abstract}
In this work, we prove the existence of a positive solution to the second-order nonlinear problem $u''+f(t,u,u')=0$ with mixed boundary conditions, where $f$ is an $L^p$-Carathéodory function satisfying certain properties. Three boundary conditions are analyzed. The proofs of the results are based on the Mawhin's coincidence degree.

\end{abstract}

\section{Introduction}

In some recent works by Guglielmo Feltrin and Fabio Zanolin, the existence of positive solutions for nonlinear problems with boundary conditions is studied. For example, in \cite{feltrin}, they prove the existence of a positive solution for the second-order nonlinear equation

$$	
u''+f(t,u,u')=0,\  0<t<T
$$
with Neumann or periodic boundary conditions. To obtain the result, they use the Mawhin coincidence degree to guarantee the existence of at least one solution. Then they apply a weak maximum principle to prove that the solution found is non-negative. Lastly, they use a strong maximum principle to prove that the obtained solution is positive.

In this work, we extend the result obtained by Guglielmo Feltrin and Fabio Zanolin to other boundary conditions, as it will be seen later. It will be apllied the same techniques used by them, namely the Mawhin coincidence degree and a maximum principle. It is worth mentioning that the maximum principle used here is slightly different from the one used by them. There is a specific section dedicated to its proof.    

Assuming that $f:[0,T]\times[0,+\infty[\,\times \,\mathbb{R}\to\mathbb{R}$ is a $L^p$-Carathéodory function, in this work, It will be studied the second-order nonlinear problem with boundary conditions

\begin{equation}
	\left\{\begin{array}{l}\label{PROBLEMA}
		u''+f(t,u,u')=0,\  0<t<T \     \\
		\mathcal{B}(u)=(0,0),
	\end{array}\right.\tag{$\cal{P}$}
\end{equation}
where the linear operator $\mathcal{B}:C^1([0,T],\R)\to\R^2$ represents the boundary conditions, which can be 

$$
\mathcal{B}(u)=\left(u'(T)-u(0),u'(0)-u(0)\right),
$$

$$
\mathcal{B}(u)=\left(u(T)-u(0),u'(0)\right)
$$

ou

$$
\mathcal{B}(u)=\left(u'(T)-u'(0),u(0)\right).
$$

	A \textit{solution} to problem \eqref{PROBLEMA} is a function $u:[0,T]\to\R$, of class $C^1$ such that $u'(t)$ is absolutely continuous and $u(t)$ satisfies \eqref{PROBLEMA} for almost every $t\in[0,T]$. We are interested in positive solutions of \eqref{PROBLEMA}, i.e., solutions $u$ such that $u(t)>0$ for all $t\in[0,T]$. However, when the problem is studied with the boundary condition $\mathcal{B}(u)=\left(u'(T)-u'(0),u(0)\right)$, we already know in advance that $u$ vanishes at $t=0$. In this case, the sought-after solution will be such that $u(t)>0$ for all $t\in\ ]0,T]$.

This work is organized as follows. In Section 2, some basic facts are recalled about Mawhin's coincidence degree. In this way, this section serves as a guide for the approach used in the study of the problems in the subsequent sections. In Section 3, It is presented the main notations used in the study of the presented problems. Sections 4, 5, and 6 are used to carefully address the problems, with each section dedicated to each boundary condition. Despite the similarity in the constructions of the problems, treating them separately and carefully allows us to better understand the differences generated by the alteration of the boundary conditions. At the end, It is included an appendix where it is presented the maximum principle, which is a fundamental result to ensure that the existing solutions are positive.

\section{Basic facts about Mawhin's coincidence degree}

Consider the real Banach spaces $E$ and $F$, and the linear Fredholm mapping $L: \operatorname{dom} L \subseteq E \rightarrow F$ of index zero. Let $\operatorname{ker} L=L^{-1}(0)$ denote the kernel or nullspace of $L$, and $\operatorname{Im} L \subseteq F$ denote the range or image of $L$. We choose linear continuous projections $P: E \rightarrow \operatorname{ker} L$ and $Q: F \rightarrow \operatorname{coker} L \subseteq F$ with $\operatorname{coker} L \cong F / \operatorname{Im} L$ being the complementary subspace of $\operatorname{Im} L$ in $F$. The linear subspace $\operatorname{ker} P \subseteq E$ is the complementary subspace of $\operatorname{ker} L$ in $E$. This gives us the direct sum decomposition 

$$E=\operatorname{ker} L \oplus \operatorname{ker} P \text{\ and\ } Z=\operatorname{Im} L \oplus \operatorname{Im} Q.$$

We define the right inverse of $L$ as $K_P: \operatorname{Im} L \rightarrow \operatorname{dom} L \cap \operatorname{ker} P$, which satisfies $L K_P(w)=w$ for each $w \in \operatorname{Im} L$. As $L$ is a Fredholm mapping of index zero, we have that $\operatorname{Im} L$ is a closed subspace of $F$, and $\operatorname{ker} L$ and $\operatorname{coker} L$ are finite dimensional vector spaces of the same dimension. Now fix an orientation on these spaces and consider a linear (orientation-preserving) isomorphism $J$ : coker $L \rightarrow \operatorname{ker} L$.

Now, let $N: E \rightarrow F$ be a possibly nonlinear operator, and consider the coincidence equation 
 
 \begin{equation}\label{intcoinc}
 	L u=N u, \quad u \in \operatorname{dom} L .
 \end{equation}
 
 According to \cite{Mawhin_livro}, the coincidence equation is equivalent to the fixed point problem 
 
 \begin{equation}\label{probpontofixo}
 	u=\Phi(u):=P u+J Q N u+K_P(I d-Q) N u, \quad u \in E,
 \end{equation}
 
 where $u \in E$. Mawhin's coincidence degree theory applies Leray-Schauder degree to the operator equation, allowing for the solution of the equation when $L$ is not invertible. To do so, it is added some structural assumptions on $N$ in order to have a completely continuous right-hand side in $\Phi$. Specifically, we assume that the operator $N$ is $L$-completely continuous, i.e., $N$ is continuous and $Q N(B)$ and $K_P(I d-Q) N(B)$ are relatively compact sets for each bounded set $B \subseteq E$. One typical scenario where the $L$-complete continuity of $N$ is satisfied is when $N$ is continuous, maps bounded sets to bounded sets, and $K_P$ is a compact linear operator.

Let $\Omega \subseteq E$ be an open and bounded set such that $L u \neq N u$ for all $u \in \operatorname{dom} L \cap \partial \Omega$. Then, the coincidence degree of $L$ and $N$ in $\Omega$ is defined as 

$$D_L(L-N, \Omega):=\operatorname{deg}(I d-\Phi, \Omega, 0),$$ 
where ``deg" denotes the Leray-Schauder degree.

Mawhin's coincidence degree theory is a powerful tool for solving equation (\ref{intcoinc}) when $L$ is not invertible. To apply this theory, It is needed to make some structural assumptions on the possibly nonlinear operator $N$. Specifically, It is assumed that $N$ is $L$-completely continuous, meaning that $N$ is continuous and that for every bounded set $B \subseteq E$, both $Q N(B)$ and $K_P(I d-Q) N(B)$ are relatively compact sets. One example of a situation where the $L$-complete continuity of $N$ is satisfied is when $N$ is continuous, maps bounded sets to bounded sets, and $K_P$ is a compact linear operator.

Suppose we have an open and bounded set $\Omega \subseteq E$ such that $L u \neq N u$ for all $u \in \operatorname{dom} L \cap \partial \Omega$. In this case, we can define the coincidence degree of $L$ and $N$ in $\Omega$ as $$D_L(L-N, \Omega) := \operatorname{deg}(I d-\Phi, \Omega, 0),$$ where ``deg" denotes the Leray-Schauder degree. Here, $\Phi$ is the operator defined in equation (\ref{probpontofixo}). We denote the (finite dimensional) Brouwer degree by $d_B$.

Remarkably, the coincidence degree of $L$ and $N$ in $\Omega$ is independent of the choice of projectors $P$ and $Q$, and is also independent of the choice of linear isomorphism $J$, provided that it is fixed an orientation on $\operatorname{ker} L$ and coker $L$ and only consider orientation-preserving isomorphisms for $J$. This generalized degree has all the usual properties of Brouwer and Leray-Schauder degree, including additivity/excision and homotopic invariance. In particular, if $D_L(L-N, \Omega) \neq 0$, then the equation (\ref{intcoinc}) has at least one solution in $\Omega$.

A crucial result for computing the coincidence degree is the relationship between the coincidence degree and the finite dimensional Brouwer degree of the operator $N$ projected into $\operatorname{ker} L$. This result was first given in \cite{Mawhin34} in its abstract form and was previously discussed in \cite{Mawhin33} in the context of periodic problems for ODEs.

\section{Notations and preliminaries}

For a fixed $T$, let $E:=C^1([0,T],\R)$ equipped with the norm
\[
\|u\|_1=\|u\|_\infty+\|u'\|_\infty
\]
and let $F:=L^1([0,T],\R)$ equipped with the norm $L_1$, denoted by $\|\cdot\|_{L_1}$.

We define $L:\dom L\to Z$ by
\[
(Lu)(t):=-u''(t),\,\,\,\,t\in[0,T],
\]
where $\dom L\subseteq E$ will be defined in each problem, since it depends on the boundary condition used.

Now we introduce a nonlinear operator $N:E\to F$, called the Nemytskii operator. First, let $f:[0,T]\times[0,+\infty[\,\times \,\mathbb{R}\to\mathbb{R}$ be an $L^p$-Carathéodory function, for some $1\leq p\leq \infty$, satisfying the following conditions:
\begin{itemize}
	\item [($f_1$)] $f(t,0,\xi)=0$, for almost every $t\in[0,T]$ and for every $\xi\in\mathbb{R}$;
	
	\item [($f_2$)] there exists a nonnegative function $k\in L^1[0,T]$ and a constant $\rho>0$ such that
	\[
	|f(t,s,\xi)|\leq k(t)(|s|+|\xi|),
	\]
	for almost every $t\in[0,T]$, for every $0\leq s\leq \rho$, and $|\xi|\leq \rho$;
	
	\item[$(f_3)$] In addition to the above assumptions, we will also suppose that $f(t,s,\xi)$ satisfies a kind of Bernstein-Nagumo condition in order to have $|u'(t)|$ bounded whenever $u(t)$ is bounded.
	
	For each $\eta>0$, there exists a continuous function
	\[
	\phi=\phi_\eta:[0,+\infty[\to [0,+\infty[, \,\,\textnormal{with}\,\,\int^{\infty}\dfrac{\xi^{\frac{p-1}{p}}}{\phi(\xi)}d\xi=\infty,
	\] 
	and a function $\psi=\psi_\eta\in L^p\left([0,T],[0,+\infty[\right)$ such that
	\[
	|f(t,s,\xi)|\leq\psi(t)\phi(|\xi|),\,\, \textnormal{for almost every}\,\,t\in[0,T],\forall\, s\in[0,\eta],\forall\,\xi\in\mathbb{R}. 
	\]
	For technical reasons, when dealing with Nagumo functions \(\phi(\xi)\) as above, we always assume that
	\[
	\liminf_{\xi\to+\infty}\phi(\xi)>0.
	\]
	This avoids the possibility of pathological examples as can be seen in \cite[p. 46-47]{Decoster} and does not affect our application.
\end{itemize}

As a first step, we extend $f$ to a Carathéodory function $\tilde{f}$ defined on $[0,T]\times \mathbb{R}^2$, by
\[
\tilde{f}(t,s,\xi)=\left\{\begin{array}{ll}
	f(t,s,\xi),&\text{if}\,\, s\geq 0\\
	-s,&\text{if}\,\, s\leq 0
\end{array}\right.
\]
and denote by $N:E\to F$ the Nemytskii operator induced by $\tilde{f}$, that is,
\[
(Nu)(t):=\tilde{f}(t,u(t),u'(t)),\,\,t\in[0,T].
\]

\section{First problem: $u'(0)=u'(T)=u(0)$}

In this section, we are interested in showing the existence of a positive solution to the problem

\begin{equation}\label{principal}
	\left\{\begin{array}{l}
		u''+f(t,u,u')=0,\,\,\,\,0<t<T\\
		u'(0)=u'(T)=u(0).		 
	\end{array}\right.
\end{equation}

A solution to problem (\ref{principal}) is a function $u:[0,T]\to\R$ of class $C^1$ such that $u'$ is absolutely continuous and $u(t)$ satisfies (\ref{principal}) for almost every $t\in[0,T]$. In this section, we are interested in solutions $u$ of (\ref{principal}) with $u(t)>0$ for all $t\in[0,T]$.

	For this problem, we define $L:\dom L\to F$ by
	\[
	(Lu)(t):=-u''(t),\,\,\,\,t\in[0,T],
	\]
	where $\dom L\subseteq E$ is the vector subspace
	\[
	\dom L=\{u\in E:u' \text{ is absolutely continuous and } u'(0)=u'(T)=u(0)\}.
	\]
	
	In these conditions, the kernel of the operator $L$ is given by
	\[
	\ker L = \{u\in E: u(t)=at+a,\, a\in \mathbb{R}\},
	\]
	which can be identified with the set of real numbers, i.e., $\ker L\equiv\R$. In fact, let $u\in \ker L$. That is, $u''(t)=0$ for all $t\in [0,T]$. Thus, $u(t)=at+b$, with $a,b\in\mathbb{R}$. Hence, $u(0)=b$ and $u'(0)=u'(T)=a$. Since $u(0)=u'(0)=u'(T)$, we have $a=b$. Therefore, $u\in \{u\in E: u(t)=at+a,\, a\in \mathbb{R}\}$, i.e., $\ker L \subseteq \{u\in E: u(t)=at+a,\, a\in \mathbb{R}\}$. On the other hand, let $u\in E$ such that $u(t)=at+a$ with $a\in\mathbb{R}$. It is clear that $u'(t)$ is absolutely continuous. Moreover, $u(0)=u'(0)=u'(T)=a$. Therefore, $u\in\dom L$. And obviously, $u''(t)=0$ for all $x\in[0,T]$. Hence, $\{u\in E: u(t)=at+a,\, a\in \mathbb{R}\}\subseteq\ker L$. We conclude that $\ker L = \{u\in E: u(t)=at+a,\, a\in \mathbb{R}\}$.
	
	Furthermore, the image of the operator $L$ is given by
	\[
	\im L=\left\{w\in F:\displaystyle\int_{0}^{T}w(t)dt=0\right\}.
	\]
	Indeed, let $w\in\im L$. Then, $w=-u''$ for some $u\in \dom L$. Thus,
	\[
	\int_{0}^{T}w(t)dt=-\int_{0}^{T}u''(t)dt=-u'(T)+u'(0)=0.
	\]
	On the other hand, let $w\in Z$ such that $\displaystyle\int_0^Tw(t)dt=0$. Define, for $s\in[0,T]$,
	\[
	v(s):=-\int_{0}^{s}w(t)dt\,\,\,\text{and}\,\,\,u(s):=\int_{0}^{s}v(t)dt.
	\]
	By the fundamental theorem of calculus, we have, for $t\in[0,T]$,
	\[
	u'(s)=v(s)\,\,\,\text{and}\,\,\,-u''(s)=-v'(s)=w(s).
	\]
	It is clear that $u'$ is absolutely continuous and, moreover,
	\[
	u'(0)=v(0)=0,\,\, u'(T)=v(T)=0\,\,\,\text{and}\,\,\, u(0)=0.
	\]
	Thus, $u\in\dom L$ and $Lu=w$. We conclude that
	\[
	\im L=\left\{w\in F:\int_{0}^{T}w(t)dt=0\right\}
	\]

	At this point, we can define the projections $P:E\to\ker L$ and $Q:F\to\operatorname{coker} L$ as follows:
	\[
	(Pu)(t)=\dfrac{u(T)-u(0)}{T}t+\dfrac{u(T)-u(0)}{T}     \textnormal{\, e\, }   (Qw)(t)=\dfrac{1}{T}\int_{0}^{T}\!\!\!w(t)dt.
	\]
	
	Since $L$ is a Fredholm operator of index zero, $\dim(\ker L) =\dim(\operatorname{coker} L) $, so $\operatorname{coker} L\equiv\mathbb{R}$. Moreover, $\ker P$ is given by $C^1$ functions such that $u(0)=u(T)$ and the linear operator $K_p:\operatorname{im} L\to\operatorname{dom} L\cap\ker P$, which is the right inverse of $L$, associates to each $w\in L^1([0,T],\R)$ with $\displaystyle	\int_{0}^{T}\!\!\!w(t)dt=0$, an unique solution $u(t)$ of
	\[
	u''+w(t)=0,\,\,\, u'(0)=u'(T)=u(0)=u(T).
	\]
	
	In this scenario, $u$ is a solution of the equation
	\begin{equation}\label{eq. conincidencia}
		Lu=Nu,\,\,u\in\dom L,
	\end{equation}
	if, and only if, it is a solution of the problem
	\begin{equation}\label{prob. equiv}
		\left\{\begin{array}{l}
			u''+\tilde{f}(t,u,u')=0,\,\,\,0<t<T\\
			u'(0)=u'(T)=u(0).
		\end{array}\right.
	\end{equation}
	Furthermore, from the definition of $\tilde{f}$ for $s<0$ and the conditions $(f_1)$ and $(f_2)$, it is easy to verify, using the maximum principle, that if $u\not\equiv 0$, then $u(t)$ is strictly positive and hence is a (positive) solution of problem $(\ref{principal})$.
	
	Next, it is presented the main result of this section, where it is proved the existence of a positive solution for the problem (\ref{principal}).
	
	\begin{theorem}\label{resultado1}
		Assume $(f_1)$, $(f_2)$, and $(f_3)$, and suppose that there exist two constants $r, R > 0$, with $r\neq R$, such that the following hypotheses are true:
		\begin{itemize}
			\item[$(H_r)$] The condition
			\[
			\int_{0}^{T}f(t,a+at,a)dt<0,\,\, \text{para}\,\,a=\dfrac{r}{1+T}
			\]
			is satisfied. Moreover, any solution $u(t)$ of the problem 
			\begin{equation}\label{probHrb1}
				\left\{\begin{array}{l}
					u''+\vartheta f(t,u,u')=0,\,\,\,0<t<T\\
					u'(0)=u'(T)=u(0),
				\end{array}\right.
			\end{equation}
			for $0<\vartheta\leq 1$, such that $u(t)>0$ in $[0,T]$, satisfies $|u|_\infty\neq r$.\\
			\item[$(H_R)$] There exist a non-negative function $v\in L^p([0,T],\R)$ with $v\not\equiv 0$ and a constant $\alpha_0>0$, such that every solution $u(t)\geq 0$ of the problem
			\begin{equation}\label{probHRb1}
				\left\{\begin{array}{l}
					u''+f(t,u,u')+\alpha v(t)=0,\,\,\,0<t<T\\
					u'(0)=u'(T)=u(0),
				\end{array}\right.
			\end{equation}
			for $\alpha\in[0,\alpha_0]$, satisfies $|u|_\infty\neq R$. Moreover, there are no solutions $u(t)$ of $(\ref{probHRb1})$ for $\alpha=\alpha_0$ with $0\leq u(t)\leq R$, for every $t\in[0,T]$.		 
		\end{itemize}
		Then the problem $(\ref{principal})$ has at least one positive solution $u(t)$ with
		\[
		\min\{r,R\}<\max_{t\in[0,T]} u(t)<\max\{r,R\}.
		\]
	\end{theorem}
	
	\begin{proof}
		As already observed, the choices of spaces $E, \,\dom L,\, F$ and operators $L:u\mapsto -u''$ and $N$ imply the equivalence between $(\ref{eq. conincidencia})$ and problem $(\ref{prob. equiv})$. It is also easy to observe that $L$ is a Fredholm operator of index zero and $N$ is $L$-completely continuous.
		
		Let us focus on the case $0<r<R$.
		
		The equation
		\begin{equation}\label{eq. conincidencia_parametrob1}
			Lu=\vartheta Nu,\,\,u\in\dom L,
		\end{equation}
		is equivalent to
		\begin{equation}\label{prob. equiv_parametrob1}
			\left\{\begin{array}{l}
				u''+\vartheta\tilde{f}(t,u,u')=0,\,\,\,0<t<T\\
				u'(0)=u'(T)=u(0).
			\end{array}\right.
		\end{equation}
		Let $u$ be a solution of $(\ref{eq. conincidencia_parametrob1})$ for some $\vartheta>0$. From the definition of $\tilde{f}$ for $s\leq 0$ and the maximum principle, we have $u(t)\geq0$ for all $t\in[0,T]$, and thus $u$ is a solution of $(\ref{probHrb1})$. Moreover, by $(f_2)$, if $u\not\equiv 0$, then $u(t)>0$ for all $t\in[0,T]$.
		
		According to condition $(f_3)$, let $\phi=\phi_r:[0,+\infty[,\to[0,+\infty[$ and $\psi=\psi_r\in L^p([0,T])$ be such that $|f(t,s,\eta)|\leq\psi(t)\phi(|\eta|)$, for almost every $t\in[0,T]$, for all $s\in[0,r]$ and $\eta\in\mathbb{R}$. By Nagumo's lemma \cite[\S\ 4.4, Proposição 4.7]{Decoster}, there is a constant $M=M_r>r$ (depending on $r$, $\phi$, and $\psi$, but not depending on $u(t)$ and $\vartheta\in, ]0,1]$) such that any solution of $(\ref{eq. conincidencia_parametrob1})$ or, equivalently, any non-negative solution of $(\ref{probHrb1})$ (for some $\vartheta\in,]0,1]$) satisfying $|u|_\infty\leq r$ is such that $|u'|_\infty\leq M_r$. Thus, condition $(H_r)$ implies that, for the open and bounded subset $\Omega_r$ of $E$ defined by
		\[
		\Omega_r:=\{u\in E:\|u\|_\infty<r,\|u'\|_\infty<M_r\},
		\]  
		we have
		\[
		Lu\neq\vartheta Nu,\,\,\forall u\in\dom L\cap \partial\Omega_r,\forall\vartheta\in\,]0,1].
		\]
		
		Let $u\in\partial\Omega_r\cap\ker L=\left\{u\in X: u(t)=a+at,\,|a|=\dfrac{r}{1+T}\right\}$. In this case,
		\[
		-JQNu=-\dfrac{1}{T}\int_{0}^{T}\tilde{f}(t,a+at,a)dt,\,\,\,\,|a|=\dfrac{r}{1+T}.
		\]
		Moreover, note that $\Omega_r\cap\ker L=\left]-\dfrac{r}{1+T},\dfrac{r}{1+T}\right[$.
		
		By the definition of $\tilde{f}$, we have
		\[
		h(a):=-\dfrac{1}{T}\displaystyle{\int}_{\!\!\!0}^{T}\tilde{f}(t,a+at,a)dt=\left\{\begin{array}{ll}
			-\dfrac{1}{T}\displaystyle{\int}_{\!\!\!0}^{\,T}f(t,a+at,a)dt,&\text{if}\,\,a>0\\
			a,&\text{if}\,\, a\leq 0
		\end{array}\right.
		\]
		Thus $QNu\neq 0$ for every $u\in\partial\Omega_r\cap\ker L$ and furthermore,
		\[
		\degb\left(h,\left]-\dfrac{r}{1+T},\dfrac{r}{1+T}\right[,0\right)=1,
		\]
		since $h\left(-\dfrac{r}{1+T}\right)<0<h\left(\dfrac{r}{1+T}\right)$. By the finite-dimensional reduction of the Mawhin coincidence degree (\cite{feltrin} - Lemma $2.1$), we conclude that
			
		\begin{equation}\label{deg1b1}
			D_L(L-N,\Omega_r)=\degb\left(h,\left]-\dfrac{r}{1+T},\dfrac{r}{1+T}\right[,0\right)=1.
		\end{equation}
		
		Now we will study the equation
		\begin{equation}\label{eqHRb1}
			Lu=Nu+\alpha v,\,\,u\in\dom L,
		\end{equation}
		for some $\alpha\geq 0$, with $v$ as in $(H_R)$. This equation is equivalent to the problem
		\begin{equation}\label{probHR1b1}
			\left\{\begin{array}{l}
				u''+\tilde{f}(t,u,u')+\alpha v(t)=0\\
				u'(0)=u'(T)=u(0).
			\end{array}\right.
		\end{equation} 
		Let $u$ be any solution of $(\ref{eqHRb1})$ for some $\alpha\geq 0$. From the definition of $\tilde{f}$ for $s\leq 0$ and the maximum principle, we have $u(t)\geq 0$ for all $t\in[0,T]$ and thus $u$ is a solution of $(\ref{probHRb1})$.
		
		According to the condition $(f_3)$, let $\phi=\phi_R:[0,+\infty[\to[0,+\infty[$ and $\psi\in L^p([0,T],\R)$ such that $|f(t,s,\xi)|\leq\psi(t)\phi(|\xi|)$, for almost every $t\in[0,T]$, for all $s\in[0,R]$ and for all $\xi\in\R$. Taking $\alpha\in[0,\alpha_0]$, we obtain
		\[
		|f(t,s,\xi)+\alpha v(t)|\leq\psi(t)\phi(|\xi|)+\alpha_0 v(t)\leq\tilde{\phi}(t)\tilde{\psi}(|\xi|)
		\] 
		for almost every $t\in[0,T]$, for all $s\in[0,R]$ and for all $\xi\in\R$, where
		\[
		\tilde{\psi}(t)=\psi(t)+\alpha_0 v(t)\,\,\,\text{and}\,\,\,\tilde{\psi}(|\xi|)=\psi(|\xi|)+1.
		\]
		Note also that \[\tilde{\psi}\in L^p([0,T]) \,\,\,\,\text{and}\,\,\,\, \int^{\infty}\dfrac{\xi^{\frac{p-1}{p}}}{\tilde{\phi}(|\xi|)}d\xi=\infty.\]
		
		Using Nagumo's lemma, there is a positive constant $M=M_R>M_r$ such that any solution of $(\ref{probHR1b1})$, or equivalently, any (nonnegative) solution of $(\ref{probHRb1})$ satisfying $|u|_\infty\leq R$ satisfies $|u'|_\infty< M_R$. Thus, the condition $(H_R)$ implies that, for the open and bounded set $\Omega_R\subseteq E$ defined by
		\[
		\Omega_R=\{u\in E:\|u\|_\infty<R, \|u'\|_\infty<M_R\}
		\]
		it holds that
		\[
		Lu\neq Nu+\alpha v,\,\,\forall u\in \dom L\cap \partial\Omega_R,\,\,\forall \alpha\in [0,\alpha_0].
		\]
		In addition, the last assumption in $(H_R)$ implies that
		\[
		Lu\neq Nu+\alpha_0 v,\,\,\forall u\in \dom L\cap \partial\Omega_R.
		\] 
		Due to the homotopy invariance of the Mawhin coincidence degree, we have
		\begin{equation}\label{deg0b1}
			D_L(L-N,\Omega_R)=0.
		\end{equation}
		
		Using $(\ref{deg1b1})$, $(\ref{deg0b1})$ and the additivity of the degree, it follows that
		\[
		D_L(L-N,\Omega_R\setminus\overline{\Omega_r})=-1.
		\] 
		This guarantees the existence of a nontrivial solution $\tilde{u}$ of $(\ref{eq. conincidencia})$ with $\tilde{u}\in \Omega_R\setminus\overline{\Omega_r}$. As $\tilde{u}$ is a nontrivial solution of $(\ref{prob. equiv})$, by the maximum principle, it follows that $\tilde{u}(t)>0$ for all $t\in [0,T]$.
		
		In the case where \(0<R<r\),  we proceed analogously. Regarding the previous case, the only relevant change is the following. First, we fix a constant \(M=M_R>0\) and obtain (\ref{deg0b1}) for the set \(\Omega_R\). Then, we repeat the first part of the proof above, fix a constant \(M_r>M_R\) and obtain (\ref{deg1b1}) for the set  \(\Omega_r\). Now we have
		\[
		D_L(L-N,\Omega_r\setminus\overline{\Omega_R})=1.
		\]
		This guarantees the existence of a nontrivial solution \(\tilde{u}\) of $(\ref{eq. conincidencia})$ whith $\tilde{u}\in\Omega_r\setminus\overline{\Omega_R}$ and we conclude, as above, that \(\tilde{u}(t)>0\) for all \(t\in[0,T]\) (by the strong maximum principle).
	\end{proof}

	\section{Second problem: $u'(0)=u'(T)$ and $u(0)=0$ }
	
	In this section, we are interested in showing the existence of a positive solution to the problem
	
	\begin{equation}\label{principal_cond2}
		\left\{\begin{array}{l}
			u''+f(t,u,u')=0,\,\,\,\,0<t<T\\
			u'(0)=u'(T)\text{ and }u(0)=0.		 
		\end{array}\right.
	\end{equation}
	
	As we will see, the construction of this section is similar to what was done in the previous section. However, the result obtained is different. Here it is not possible to obtain a solution that is positive for all $t\in[0,T]$, since we have $u(t)=0$ in the boundary condition. However, this will be the only point where the solution vanishes in the interval $[0,T]$.
	
	A solution to problem (\ref{principal_cond2}) is a function $u:[0,T]\to\R$ of class $C^1$ such that $u'$ is absolutely continuous and $u(t)$ satisfies (\ref{principal}) for almost every $t\in[0,T]$. As already noted, in this section, we are interested in solutions $u$ of (\ref{principal_cond2}) with $u(t)>0$ for all $t\in\ ]0,T]$.

	For this problem, we define $L:\dom L\to F$ by
	\[
	(Lu)(t):=-u''(t),\,\,\,\,t\in[0,T],
	\]
	where $\dom L\subseteq E$ is the vector subspace
	\[
	\dom L=\{u\in E:u' \text{ is absolutely continuous and }  u'(0)=u'(T)\text{ and }u(0)=0\}.
	\]
	
	In these conditions, the kernel of the operator $L$ is given by
	\[
	\ker L = \{u\in E: u(t)=at,\, a\in \mathbb{R}\},
	\]
	which can be identified with the set of real numbers, i.e, \(\ker L\equiv\R\). In fact,  let $u\in \ker L$. That is, $u''(t)=0$ for all $t\in [0,T]$. Thus, $u(t)=at+b$ with $a,b\in\mathbb{R}$. Therefore, $u(0)=b$ and $u'(0)=u'(T)=a$. Since $u(0)=0$, we have $b=0$. Hence, $u\in \{u\in E: u(t)=at,, a\in \mathbb{R}\}$, i.e., $\ker L \subseteq \{u\in E: u(t)=at, a\in \mathbb{R}\}$. On the other hand, let $u\in E$ such that $u(t)=at$ with $a\in\mathbb{R}$. It is clear that $u'(t)$ is absolutely continuous. Moreover, $u(0)=0$ and $u'(0)=u'(T)=a$. Therefore, $u\in\dom L$. And obviously, $u''(t)=0$ for all $t\in[0,T]$. Hence, $\{u\in E: u(t)=at,, a\in \mathbb{R}\}\subseteq\ker L$. We conclude that $\ker L = \{u\in E: u(t)=at,, a\in \mathbb{R}\}$.
	Furthermore, the image set of the operator $L$ is given by
	\[
	\im L=\left\{w\in F:\displaystyle\int_{0}^{T}w(t)dt=0\right\}.
	\]
	In fact, let $w\in\im L$. Then, $w=-u''$ for some $u\in \dom L$. Thus, 
	\[
	\int_{0}^{T}w(t)dx=-\int_{0}^{T}u''(t)dt=-u'(T)+u'(0)=0.
	\]
	On the other hand, let $w\in F$ such that $\displaystyle\int_0^Tw(t)dt=0$. Define, for $s\in[0,T]$,
	\[
	v(s):=-\int_{0}^{s}w(t)dt\,\,\,\text{e}\,\,\,u(s):=\int_{0}^{s}v(t)dx.
	\]
	
	By the fundamental theorem of calculus, we have, for $s\in[0,T]$,
	\[
	u'(s)=v(s)\,\,\,\text{e}\,\,\,-u''(s)=-v'(s)=w(s).
	\]
	Clearly, $u'$ is absolutely continuous and, furthermore,
	\[
	u'(0)=v(0)=0,\,\, u'(T)=v(T)=0\,\,\,\text{e}\,\,\, u(0)=0.
	\]
	Thus, $u\in\dom L$ and $Lu=w$. We conclude that
	\[
	\im L=\left\{w\in F:\int_{0}^{T}w(t)dt=0\right\}.
	\]
	
	At this point, we can define the projections $P: E\to\ker L$ and $Q: F\to\coker L$ as follows: 
	
	\[
	(Pu)(t)=u'(0)t     \textnormal{\, and\, }   (Qw)(t)=\dfrac{1}{T}\int_{0}^{T}\!\!\!w(t)dt.
	\]
	
	As $L$ is a Fredholm operator of index zero, $\dim(\ker L) =\dim(\coker L) $, so that $\coker L\equiv\mathbb{R}$. Moreover, $\ker P$ consists of $C^1$ functions such that $u'(0)=0$, and the linear operator $K_p:\im L\to\dom\cap\ker P$ associates to each $w\in L^1([0,T],\R)$ with $\displaystyle\int_{0}^{T}\!\!\!w(t)dt=0$, a unique solution $u(t)$ of
	\[
	u''+w(t)=0,\,\,\, u'(0)=u'(T)=u(0)=0.
	\]
	
	In this scenario, $u$ is a solution of the equation
	\begin{equation}\label{eq. conincidenciaprob2}
		Lu=Nu,\,\,u\in\dom L,
	\end{equation}
	if, and only if, it is a solution of the problem
	\begin{equation}\label{prob2. equiv}
		\left\{\begin{array}{l}
			u''+\tilde{f}(t,u,u')=0,\,\,\,0<t<T\\
			u'(0)=u'(T)\text{ and }u(0)=0.
		\end{array}\right.
	\end{equation}
	In addition, from the definition of $\tilde{f}$ for $s<0$ and conditions $(f_1)$ and $(f_2)$, it is easy to verify, using the maximum principle, that if $u\not\equiv 0$, then $u(t)$ is strictly positive and thus is a (positive) solution of problem $(\ref{principal_cond2})$.
		
	Now we present the main result of this section, which guarantees the existence of a positive solution to problem (\ref{principal_cond2}). The proof of this result will be omitted, as it is analogous to that of Theorem \ref{resultado1}. The difference we can highlight, besides the boundary condition, lies in the hypothesis $(H_r)$ (integral condition). This difference is related to the kernel of the operator $L$.
	
	\begin{theorem}
		Assume $(f_1)$, $(f_2)$, and $(f_3)$, and suppose that there exist two constants $r,R>0$, with $r\neq R$, such that the following hypotheses are true.	
		\begin{itemize}
			\item[$(H_r)$] The condition
			\[
			\int_{0}^{T}f\left(t,at,a\right)dt<0,\,\, \text{para}\,\,a=\dfrac{r}{T}
			\]
			is satisfied. Moreover, any solution $u(t)$ of the problem
			\begin{equation}\label{probHrb2}
				\left\{\begin{array}{l}
					u''+\vartheta f(t,u,u')=0,\,\,\,0<t<T\\
					u'(0)=u'(T)\text{ and }u(0)=0,
				\end{array}\right.
			\end{equation}
			for $0<\vartheta\leq 1$, such that $u(t)>0$ in $[0,T]$, satisfies $\|u\|_\infty\neq r$.\\
			\item[$(H_R)$] here exists a non-negative function $v\in L^p([0,T],\R)$ with $v\not\equiv 0$ and a constant $\alpha_0>0$, such that every solution $u(t)\geq 0$  of the problem
			\begin{equation}\label{probHRb2}
				\left\{\begin{array}{l}
					u''+f(t,u,u')+\alpha v(t)=0,\,\,\,0<t<T\\
					u'(0)=u'(T)\text{ and }u(0)=0,
				\end{array}\right.
			\end{equation}
			for $\alpha\in[0,\alpha_0]$, satisfies $\|u\|_\infty\neq R$. Moreover, there are no solutions $u(t)$ of $(\ref{probHRb2})$ for $\alpha=\alpha_0$ with $0\leq u(t)\leq R$, for all $t\in[0,T]$.		 
		\end{itemize}
		Then the problem $(\ref{principal_cond2})$  has at least one positive solution $u(x)$ with
		\[
		\min\{r,R\}<\max_{t\in[0,T]} u(t)<\max\{r,R\}.
		\]
	\end{theorem}

	\section{Third problem: $u(0)=u(T)\text{ and }u'(0)=0$ }
	
	In this section, we are interested in showing the existence of a positive solution to the problem
	
	\begin{equation}\label{principal_cond3}
		\left\{\begin{array}{l}
			u''+f(t,u,u')=0,\,\,\,\,0<t<T\\
			u(0)=u(T)\text{ and }u'(0)=0.		 
		\end{array}\right.
	\end{equation}
	
	Here we will again make a construction similar to what was seen in previous sections. However, with the boundary conditions we will use, we will see a significant difference in the range of the operator $L$ and, consequently, in the projections $P$ and $Q$.
	
	A solution to problem (\ref{principal_cond3}) is a function \(u:[0,T]\to\R\) of class \(C^1\) such that \(u'\) is absolutely continuous and \(u(t)\) satisfies (\ref{principal_cond3}) for almost every $t\in[0,T]$. In this section, we are interested in solutions $u$ of (\ref{principal_cond3})
	with $u(t)>0$ for all \(t\in\ [0,T]\).
	
	For this problem, we define $L:\dom L\to F$ by
	\[
	(Lu)(t):=-u''(t),\,\,\,\,t\in[0,T],
	\]
	where $\dom L\subseteq E$ is the vector subspace
	\[
	\dom L=\{u\in E:u' \text{ is absolutely continuous and } u(0)=u(T)\text{ and }u'(0)=0\}.
	\]
	
	In these conditions, the kernel of the operator $L$ is given by
	\[
	\ker L = \{u\in E: u(t)=b,\, b\in \mathbb{R}\},
	\]
	which can be identified with the set of real numbers, i.e., $\ker L\equiv\R$. In fact, let $u\in \ker L$. $u''(t)=0$ for all $t\in [0,T]$. Thus, $u(t)=at+b$, with $a,b\in\mathbb{R}$. Hence, $u(0)=b$ e $u'(0)=a$. Since $u'(0)=0$, we have $a=0$. Therefore, $u\in \{u\in E: u(t)=b,\, b\in \mathbb{R}\}$, i.e, $\ker L \subseteq \{u\in E: u(t)=b,\, b\in \mathbb{R}\}$. On the other hand, let $u\in E$ such that $u(t)=b$ with $b\in\mathbb{R}$. It is clear that $u'(t)$ is absolutely continuous. Moreover, $u(0)=u(T)=b$ and $u'(0)=0$. Portanto, $u\in\dom L$. Therefore, $u\in\dom L$. And obviously, $u''(t)=0$ for all $t\in[0,T]$. Hence, $\{u\in E: u(t)=b,\, b\in \mathbb{R}\}\subseteq\ker L$. We conclude that $\ker L = \{u\in X: u(t)=b,\, b\in \mathbb{R}\}$. 
	  
	Furthermore, the image of the operator $L$ is given by
	\[
	\im L=\left\{w\in F:\displaystyle\int_{0}^{T}\left(\int_{0}^{s}w(t)dt\right)ds=0\right\}
	.\]
	Indeed, let $w\in\im L$. Then, $w=-u''$ for some $u\in \dom L$. Thus, 
	\[
	\int_{0}^{T}\left(\int_{0}^{s}w(t)dt\right)ds=-\int_{0}^{T}\left(\int_{0}^{s}u''(t)dt\right)ds=-\int_{0}^{T}u'(s)ds=-u(T)+u(0)=0.
	\]
	On the other hand, let $w\in F$ such that $\displaystyle\int_0^Tw(t)dt=0$. Define, for $t\in[0,T]$,
	\[
	u(t):=u(0)-\int_{0}^{t}\left(\int_{0}^{x}w(s)ds\right)dx.
	\].
	
	By the fundamental theorem of calculus, we have, for $t\in[0,T]$,
	\[
	u'(t)=-\int_{0}^{t}w(s)ds\,\,\,\text{and}\,\,\,-u''(t)=w(t).
	\]
	It is clear that $u'$ is absolutely continuous and, moreover,
	\[
	u'(0)=0\text{ and } u(0)=u(T).
	\]
	Thus, $u\in\dom L$ and $Lu=w$. We conclude that
	\[
	\im L=\left\{w\in F:\int_{0}^{T}\left(\int_{0}^{s}w(t)dt\right)ds=0\right\}.
	\].
	
	At this point we can define the projections $P: E\to\ker L$ e $Q: F\to\coker L$ as follows: 
	
	\[
	(Pu)(t)=\dfrac{1}{T}\int_{0}^{T}u(t)dt     \textnormal{\, and\, }   (Qw)(t)=\dfrac{2}{T^2}\int_{0}^{T}\left(\int_{0}^{s}w(t)dt\right)ds.
	\]
	
	Since $L$ is a Fredholm operator of index zero, $\dim(\ker L) =\dim(\operatorname{coker} L) $, so $\operatorname{coker} L\equiv\mathbb{R}$. Moreover, $\ker P$ is given by $C^1$ functions with mean value zero and the linear operator $K_p:\im L\to\dom\cap\ker P$, which is the right inverse of $L$, associates to each $w\in L^1([0,T],\R)$ with $\displaystyle\int_{0}^{T}\left(\int_{0}^{s}w(t)dt\right)ds=0$, a unique solution $u(t)$ of
	\[
	u''+w(t)=0,\,\,\, u(0)=u(T),\,\,\, u'(0)=0 \textnormal{ and } \dfrac{1}{T}\int_{0}^{T}u(t)dt.
	\]
	
	In this scenario, $u$ is a solution of the equation
	\begin{equation}\label{eq. conincidenciaprob3}
		Lu=Nu,\,\,u\in\dom L,
	\end{equation}
	if, and only if, it is a solution of the problem
	\begin{equation}\label{prob3. equiv}
		\left\{\begin{array}{l}
			u''+\tilde{f}(t,u,u')=0,\,\,\,0<t<T\\
			u(0)=u(T)\text{ and }u'(0)=0.
		\end{array}\right.
	\end{equation}
	Furthermore, from the definition of $\tilde{f}$ for $s<0$ and the conditions $(f_1)$ and $(f_2)$, it is easy to verify, using the maximum principle, that if $u\not\equiv 0$, then $u(t)$ is strictly positive and hence is a (positive) solution of problem (\ref{principal_cond3}).
	
Now we present the main result of this section, which guarantees the existence of a positive solution to problem (\ref{principal_cond3}). Again, there is no need to present the proof, as it is similar to what was done in Theorem \ref{resultado1}. However, the result is interesting because we can obtain positive solutions over the entire interval $[0,T]$, just like in problem $(\ref{principal})$, but with a construction that has some significant differences, as we have seen in this section.

	\begin{theorem}
		Assume $(f_1)$, $(f_2)$, and $(f_3)$, and suppose that there exist two constants $r, R > 0$, with $r\neq R$, such that the following hypotheses are true:	
		\begin{itemize}
			\item[$(H_r)$] The condition
			\[
			\displaystyle\int_{0}^{T}\left(\int_{0}^{s}f(t,r,0)dt\right)ds<0
			\]
		is satisfied. Moreover, any solution $u(t)$ of the problem
			\begin{equation}\label{probHrb3}
				\left\{\begin{array}{l}
					u''+\vartheta f(t,u,u')=0,\,\,\,0<x<T\\
					u(0)=u(T)\text{ and }u'(0)=0,
				\end{array}\right.
			\end{equation}
			for $0<\vartheta\leq 1$, such that $u(t)>0$ in $[0,T]$, satisfies $\|u\|_\infty\neq r$.\\
			\item[$(H_R)$] There exist a non-negative function $v\in L^p([0,T],\R)$ with $v\not\equiv 0$ and a constant $\alpha_0>0$, such that every solution $u(t)\geq 0$  of the problem
			\begin{equation}\label{probHRb3}
				\left\{\begin{array}{l}
					u''+f(t,u,u')+\alpha v(t)=0,\,\,\,0<t<T\\
					u(0)=u(T)\text{ and }u'(0)=0,
				\end{array}\right.
			\end{equation}
			for $\alpha\in[0,\alpha_0]$, satisfies $\|u\|_\infty\neq R$. Moreover, there are no solutions $u(t)$ of $(\ref{probHRb3})$ for $\alpha=\alpha_0$ with $0\leq u(t)\leq R$, for all $t\in[0,T]$.		 
		\end{itemize}
		Then the problem $(\ref{principal_cond3})$ has at least one positive solution $u(t)$ with
		\[
		\min\{r,R\}<\max_{t\in[0,T]} u(t)<\max\{r,R\}.
		\]
	\end{theorem}

	\section*{Appendix. Maximum principle}
	
	Consider the problem
	
	\begin{equation}\label{Eq_princ_max}
		\left\{\begin{array}{l}
			u''+h(t,u)=0\\
			\mathcal{B}(u)=(0,0),
		\end{array}\right.
	\end{equation}
	where 
	\begin{equation}\label{bordo}
		\begin{array}{ll}
			\mathcal{B}(u)\in&\!\!\!\!\Big\{\left(u'(T)-u(0),u'(0)-u(0)\right),\left(u(T)-u(0),u'(0)\right),\\ &\left(u'(T)-u'(0),u(0)\right)\Big\}.		 
		\end{array}
	\end{equation}

	\begin{lema}[Maximum principle]
		Let  $h:[0,T]\times \R\to\R$ be an  $L^1$ - Carathéodory function.
		\begin{enumerate}[(i)]
			\item If $h(t,s)>0$, almost every $t\in[0,T]$ and for all $s<0$,
			then any solution of $(\ref{Eq_princ_max})$  is non-negative on $[0,T]$.\\
			\item If $h(t,0)\equiv 0$ and there exists $q\in L^1\left([0,T],[0,+\infty[\right)$ such that
			\[
			\limsup_{s\to 0^+}\dfrac{|h(t,s)|}{s}\leq q(t),
			\]
			uniformly for almost every $t\in [0,T]$, then any non-trivial non-negative solution of $(\ref{Eq_princ_max})$ satisfies:
			\begin{itemize}
				\item $u(t)>0$ for every $t\in\,\, ]0,T]$ if $\mathcal{B}(u)=\left(u'(T)-u'(0),u(0)\right)$;
				\item $u(t)>0$ for every $t\in [0,T]$ if $\mathcal{B}(u)=\left(u(T)-u(0),u'(0)\right)$;
				\item $u(t)>0$ for every $t\in [0,T]$ if $\mathcal{B}(u)=\left(u'(T)-u(0),u'(0)-u(0)\right)$.
			\end{itemize}  
			
		\end{enumerate}	
	\end{lema}
	\begin{proof}
		\begin{enumerate}[(i)]
			\item Suppose there is a solution $u(t)$ of (\ref{Eq_princ_max}) and $\hat{t}\in[0,T]$ such that $u(\hat{t})<0$. Let $]t_1,t_2[\ \subseteq\ ]0,T[$  be the maximal
			open interval containing $\hat{t}$ such that $u(t)<0$ for all $t\in\ ]t_1,t_2[$. We will analyze four cases, namely:
			\begin{enumerate}[(a)]
				\item $0<t_1<t_2<T$:\\
				By hypothesis, $u''(t)=-h(t,u(t))<0$ for almost every $t\in\ ]t_1,t_2[$. This means that $u(t)$ is a concave function on $[t_1,t_2]$. Moreover,  $u(t_1)=u(t_2)=0$. Therefore, we have $u(t)\geq 0$ for all $t\in\ ]t_1,t_2[$,  which contradicts the assumption that $u(\hat{t})<0$. Thus, this case is not possible. 
				
				\item $t_1=0$ and $t_2=T$:\\
				By hypothesis, we know that $u''(t)=-h(t,u(t))<0$ for almost every $t\in[0,T]$. This guarantees that $u'(t)$ is strictly decreasing in $[0,T]$. Thus, we have $u'(T)<u'(0)$ which contradicts the boundary conditions in the cases $\mathcal{B}(u)=\left(u'(T)-u(0),u'(0)-u(0)\right)$ and $\mathcal{B}(u)=\left(u'(T)-u'(0),u(0)\right)$. 
				
				For the case where $\mathcal{B}(u)=\left(u(T)-u(0),u'(0)\right)$, we use the facts that $u'(t)$ is strictly decreasing in $[0,T]$ and that $u'(0)=0$ to conclude that $u'(t)<0$ for all $t\in\ ]0,T]$. Therefore, $u(t)$ is decreasing in $[0,T]$ implying that $u(T)<u(0)$, which contradicts the boundary condition. Hence, this case is also not possible.

				\item $t_1=0$ and $t_2<T$:\\
				
				By hypothesis, we know that $u''(t)=-h(t,u(t))<0$ for almost all $t\in[0,t_2]$. Therefore,
				\begin{equation}
					0>\int_{0}^{t_2}u''(t)dt=u'(t_2)-u'(0).
				\end{equation}
				Thus, we have $u'(t_2)<u'(0)$. When the boundary condition is $\mathcal{B}(u)=\left(u'(T)-u(0),u'(0)-u(0)\right)$, we have $u'(t_2)<u'(0)=u(0)\leq 0$. This is a contradiction since $u'(t_2)\geq 0$. In the case where the boundary condition is $\mathcal{B}(u)=\left(u(T)-u(0),u'(0)\right)$, we have $u'(t_2)<u'(0)=0$, again contradicting the fact that $u'(t_2)\geq 0$. We still need to look at the condition $\mathcal{B}(u)=\left(u'(T)-u'(0),u(0)\right)$. But in this case, $u(0)=u(t_2)=0$, which can be solved with the same argument used in case (a).

				\item $0<t_1$ and $t_2=T$:\\
				We start with the condition $\mathcal{B}(u)=\left(u(T)-u(0),u'(0))\right)$. Note that $u(T)\leq 0$. The possibility $u(T)=0$ can be excluded by the same strategy used in part (a). Now, if $u(T)$ is negative, by the boundary condition, $u(0)$ is also negative. Suppose $\overline{t}_2=\sup\left\{t:u\textnormal{ is negative in } [0,t[\right\}$ and use the same argument as in part (c) for the interval $[0,\overline{t}_2]$.
				
				We move to the condition $\mathcal{B}(u)=\left(u'(T)-u'(0),u(0)\right)$. Note that $u(t_1)\leq 0$. First, assume that $u'(T)\geq 0$. Since $u''(t)=-h(t,u(t))<0$ for almost every $t\in[0,t_2]$, it follows that
				
				\begin{equation}\label{bordoint}
					0>\int_{t_1}^{T}u''(t)dt=u'(T)-u'(t_1),
				\end{equation}
				which implies $u'(t_1)>u'(T)\geq 0$. This leads to a contradiction. If $u'(T)$ is negative, then by the boundary condition, we have $u(0)=0$ and $u'(0)<0$. Therefore, there is $s\in\,\,]0,T[$ such that $u(t)<0$ in $]0,s[$. Suppose \[\overline{s}=\sup\left\{s:u(t)\textnormal{ é negativa em } ]0,s[\right\}.\] Note that $\overline{s}\leq t_1<T$, so we go back to case (c).
				
				Finally, let us consider the condition $\mathcal{B}(u)=\left(u'(T)-u(0),u'(0)-u(0))\right)$. Assuming $u'(T)\geq 0$ and using $(\ref{bordoint})$ we conclude that $u'(t_1)>u'(T)\geq 0$, which leads to a contradiction. Now, if  $u'(T)<0$, then by the boundary condition, we have $u'(T)=u'(0)=u(0)<0$. Therefore, there exists $s\in\,\,]0,T[$ such that $u(t)<0$ in $]0,s[$. Suppose \[\overline{s}=\sup\left\{s:u(t)\textnormal{ is negative in } ]0,s[\right\}.\]  Note that  $\overline{s}\leq t_1<T$, so we go back to case (c).

			\end{enumerate}
			\item By contradiction, suppose that there exists a solution $u(t)\geq 0$ of  (\ref{Eq_princ_max}) and $t^*\in[0,T]$ such that $u(t^*)=0$. It is clear that if $t^*\in\ ]0,T[$, then $u'(t^*)=0$. Let us examine what happens to $u'(t^*)$ in the cases $t^*=0$ or $t^*=T$ for each boundary condition:
			\begin{itemize}
				\item $\mathcal{B}(u)=\left(u'(T)-u'(0),u(0)\right)$: In this case, we cannot make any assertions about $t^*=0$. But if $t^*=T$, since $u'(0)=u'(T)$, $u(T)=0$ and $u(t)\geq 0$ in $[0,T]$,  we conclude that  $u'(T)=0$.
				\item $\mathcal{B}(u)=\left(u(T)-u(0),u'(0)\right)$: Here, $u'(0)=0$. Therefore, we know what happens when $t^*=0$. But we cannot make any assertions for $t^*=T$.
				\item $\mathcal{B}(u)=\left(u'(T)-u(0),u'(0)-u(0)\right)$: Assuming $t^*=0$,  by the boundary conditions, we have $u'(0)=u(0)=0$. Setting $t^*=T$, then $u(T)=0$. Since $u(t)\geq 0$ for all $t\in[0,T]$, we have $u'(T)\leq 0$. By the boundary condition, $u(0)=u'(T)\leq 0$. But $u(t)\geq 0$ in $t\in[0,T]$, so $u(0)=0$, i.e., $u'(T)=0$. 
			\end{itemize}
			
		\end{enumerate}
		
		Now we will prove the following assertion:
		\begin{af}\label{Af1}
			There exist $\varepsilon>0$ such that $u(t)=0$ for all $x\in[t^*-\varepsilon,t^*+\varepsilon]$.
		\end{af}

		From the hypothesis, we know that there is $\delta>0$ such that  
		\[
		|h(t,s)|\leq q_1(t)s,\textnormal{ for almost every } t\in[0,T], \textnormal{ for all } s\in[0,\delta],
		\]
		where $q_1(t):=q(t)+1$. Using the continuity of $u(t)$, we fix $\varepsilon>0$ such that $0\leq u(t)\leq \delta$, for all $t\in [t^*-\varepsilon,t^*+\varepsilon]$. We will use $\|(\xi_1,\xi_2)\|=|\xi_1|+|\xi_2|$ as the norm for $\mathbb{R}^2$. For all $t\in\,\,]t^*,t^*+\varepsilon]$ we have	
		\begin{eqnarray*}
			0\leq\|(u(t),u'(t))\|&=&|u(t)|+|u'(t)|\\
			& = &u(t^*)+{\displaystyle \int_{t^*}^{t}u'(\xi)d\xi}+\left|u'(t^*)+{\displaystyle \int_{t^*}^{t}-h(\xi,u(\xi))d\xi}\right|\\
			& \leq & {\displaystyle \int_{t^*}^{t}|u'(\xi)|d\xi}+{\displaystyle \int_{t^*}^{t}|h(\xi,u(\xi))|d\xi}\\
			& \leq & {\displaystyle \int_{t^*}^{t}\left[q_1(\xi)|u(\xi)|+|u'(\xi)|\right]d\xi}\\
			& \leq & {\displaystyle \int_{t^*}^{t}(q_1(\xi)+1)(|u(\xi)|+|u'(\xi)|)d\xi.} 
		\end{eqnarray*}
		Using the Gronwall inequality, we obtain
		\[
		0\leq u(t)\leq \|(u(t),u'(t))\|=0,\,\,\forall t\in\,\,]t^*,t^*+\varepsilon].
		\]
		With a similar calculation, we can prove $u(t)=0$ for all $t\in[t^*-\varepsilon,t^*[$. Thus, the assertion is proved.
		
		Finally, we need to show that Assertion \ref{Af1} leads to a contradiction. We will see this in the following assertion:
		
		\begin{af}
			The Assertion \ref{Af1} implies that $u\equiv 0$ in $[0,T]$, which is a contradiction. Indeed, suppose $t^*\in\,\,]0,T[$. Seja $\varepsilon>0$ such that $J:=[t^*-\varepsilon,t^*+\varepsilon]\subseteq\,\, ]0,T[$ e $u\equiv 0$ in $J$. Consider $E=\{t\in[0,T]: u\equiv 0 \textnormal{ in } [t^*,t] \textnormal{ or } u\equiv 0 \textnormal{ in } [t,t^*]\}$.  We want to show that $E=[0,T]$. If $E\subsetneq[0,T]$, then $\sup E<T$ or  $\inf E >0$. Let $b=\sup E$ and assume $b<T$. I claim that $b\in E$. In fact, since $b=\sup E$, there is a sequence $(t_n)_n\subset E$ whith $t_n\to b$. As $u$ is continuous, ${\displaystyle\lim_{n\to\infty}u(t_n)}=u(b)=0$, since $u(t_n)=0$ for all $n$. Note that $b\in E$, otherwise there would exist $y\in\,\,]t^*,b[$ with $u(y)\neq 0$. However, as $t_n\to b$, there is $n$ such that $t_n>y$, which is absurd, since $u\equiv 0$ in $[t^*,t_n]$. Now we know that $b\in]t^*,T[$. By Assertion\ref{Af1}, there exists $\varepsilon>0$ such that $u$ is zero in $[b-\varepsilon,b+\varepsilon]$. But then $[t^*,b+\varepsilon]\subseteq E$, which is absurd, since $b=\sup E$. Hence $b=T$.
		\end{af}	
	\end{proof}
	
	\begin{obs}
		The maximum principle presented also holds for the more general problem
		\[
		\left\{\begin{array}{l}
			u''+\tilde{f}(t,u,u')=0,\ 0<t<T,\\
			\mathcal{B}(u)=0.
		\end{array}\right.
		\]
		where \(\tilde{f}:[0,T]\times\R\times\R\to\R\) is a $L^P$-Carathéodory function as in Section \(2\), taking  \(-s\) for \(s\leq 0\) and satisfying the conditions \((f_1)\) e \((f_2)\). The proof of this result is the same as seen above with small changes.
	\end{obs}

\end{document}